\documentclass[12pt]{article} \textheight 8in\textwidth
6in \oddsidemargin 0in\evensidemargin 0in
\usepackage{graphicx}
\usepackage{amsfonts}
\usepackage{epsf}

\newtheorem{theorem}{Theorem}

\newtheorem{conjecture}[theorem]{Conjecture}

\newtheorem{problem}[theorem]{Problem}

\def\beq{\begin{equation}}\def\eeq{\end{equation}}
\def\beqn{\begin{eqnarray}}\def\eeqn{\end{eqnarray}}
\def\pont{\hspace{-6pt}{\bf.\ }}

\def\qed{\ifhmode\unskip\nobreak\fi\quad\ifmmode\Box\else$\Box$\fi}

\newcommand{\text}[1]{\quad\mbox{#1}\quad}

\newcommand{\comment}[1]{}

\def\beq{\begin{equation}}\def\eeq{\end{equation}}
\def\beqn{\begin{eqnarray}}\def\eeqn{\end{eqnarray}}
\def\pont{\hspace{-6pt}{\bf.\ }}

\def\qed{\ifhmode\unskip\nobreak\fi\quad\ifmmode\Box\else$\Box$\fi}

\title{Problems and Memories\thanks{A talk at Erd\H os 100}}

\author{Andr\'as  Gy\'arf\'as \thanks{Supported in part by OTKA grant K104343}\\
Alfr\'ed R\'enyi Institute of Mathematics \\[-0.8ex]}
\begin{document}
\maketitle

\section{First encounter}

My first encounter with Paul Erd\H os took place in 1962 at the M\'atrah\'aza Guest House of the Hungarian Academy of Sciences. I was a high school student and was rather embarrassed by the solemn formalities, especially at the dinner tables. I was pleased to hear the signs of some unaccepted behavior from a neighboring table, where an ``old'' man about fifty was regulated by his mother: ``Pali you should keep your fork properly!'' I soon learned that the unruled boy is a famous mathematician who travels around the world with his mother.  Next day I had opportunity to play ping pong against ``Pali'' and I was very angry upon being beaten by such an old man playing in a ridiculous style. As a consolation he told me what a graph is and what does Tur\'an theorem say about the number of edges in a graph that does not contain $K_{k+1}$. ``Adding any edge to the Tur\'an graph we get a $K_{k+1}$ but what is the smallest graph with this property?'' -   asked during the revenge game which I lost again. Then he left but the problem bugged me and about a year later, after  a lecture he gave in Budapest, I handed him my typewritten solution for the case when the number of vertices is large in terms of $k$. It was disappointing to learn that his question was not an unsolved problem, but his result with Hajnal and Moon, \cite{EHM}. Two years later, at E\"otv\"os University, I listened to B\'ela Bollob\'as' talk at the Haj\'os seminar about extending the result to hypergraphs. I remember B\'ela's comment: ``this is trivial'' \cite{B1}...  Trivial or not, certainly important and rediscovered several times \cite{JP,KA}.  The underlying idea, cross-intersecting sets, developed further and have many applications.

\section{Memphis Tennessee...}

\noindent Long distance information, give me Memphis Tennessee,

\noindent Help me find the party trying to get in touch with me,

\noindent She could not leave her number, but I know who placed the call,

\noindent 'Cause my uncle took the message and he wrote it on the wall...

... sang Chuck Berry, and beside Elvis Presley, Chuck Berry and Fedex, Memphis Tennessee is also known as a hub during the movements of Paul Erd\H os in the US. Indeed, I have had much more chance to meet (and think on math) with him there than in Budapest. The University of Memphis (with Faudree, Ordman, Rousseau and Schelp as permanent faculty and me as permanent visitor) provided many opportunities to pursue theorems, problems and conjectures. From the early 90-s the department was fortified by B\'ela Bollob\'as, who leads a Chair of Excellence in Combinatorics since then.

\subsection{Cyles in graphs without proper subgraphs of minimum degree $3$}
The following observation is due to Erd\H os, Faudree, Rousseau and Schelp \cite{EFRS}.
\begin{itemize}\item
Graphs with $n$ vertices and $2n-1$ edges must contain proper subgraphs of minimum degree $3$ but this fails for graphs with $n$ vertices and $2n-2$ edges, for example the wheel is such a graph.
\end{itemize}

In \cite{EFGS2} we had a closer look on the family $G(n)$ of graphs with $n$ vertices, $2n-2$ edges and without proper subgraphs of minimum degree $3$. We showed that graphs in $G(n)$ contain cycles $C_3,C_4,C_5$ and also $C_k$ for $k\ge \log_2{n}$ but not necessarily for $k\ge c\sqrt{n}$. However, we could not resolve the following.

\begin{conjecture}\pont (\cite{EFGS2}) Every $G\in G(n)$ contains cycles of length $i$ for every integer $3\le i \le k$ where $k$ tends to infinity with $n$.
\end{conjecture}

\subsection{Large chordal subgraphs} In \cite{EGYOZ} we discussed the size of chordal (interval, threshold) subgraphs present in graphs with several ranges of edges. One particular jump is observed at graphs with $n$ vertices and $n^2/3$ edges.
\begin{itemize} \item Any graph with $n$ vertices and at least $n^2/3$ edges contains a chordal subgraph with at least $2n-3$ edges. The complete tripartite graph shows that this is sharp.
\end{itemize}

\begin{conjecture}\pont(\cite{EGYOZ}) Any graph with $n$ vertices and more than $n^2/3$ edges contains a chordal subgraph with at least $8n/3-4$ edges. The complete tripartite graph with one additional edge shows that this would be sharp.
\end{conjecture}

We could prove a weaker result,  that graphs  with $n$ vertices and more than $n^2/3$ edges contain chordal subgraphs with at least $7n/3-6$ edges \cite{EGYOZ}.

\subsection{Monochromatic domination} The following result, conjectured by Erd\H os and Hajnal \cite{EH}, was proven in \cite{EFGS} and independently by Kostochka \cite{KO}. In an edge colored complete graph $K$ a subset $S$ of vertices dominate in color $i$ those vertices in $V(K)-S$ that send at least one edge of color $i$ to $S$.

\begin{itemize} \item (\cite{EFGS}) Assume that the edges of $K_n$ are $2$-colored and $k$ is a positive integer. There exist $k$ vertices of $K_n$ such that in one of the colors they dominate all but at most ${n-1\over 2^k}$ vertices of $K_n$.
\end{itemize}

 It turned out, that for $3$ colors the situation is different.

\begin{itemize} \item (\cite{EFGGRS}) If the edges of $K_n$ are $3$-colored then in one of the colors at most $22$ vertices dominate at least $2n/3$ vertices of $K_n$. 
\end{itemize}

\begin{conjecture}\pont\label{threecol}(\cite{EFGGRS}) If the edges of $K_n$ are $3$-colored then in one of the colors at most $3$ vertices dominate at least $2n/3$ vertices of $K_n$. 
\end{conjecture}

The random $3$-coloring shows that no two vertices dominate much more than $5n/8$ vertices.
Note that one cannot dominate more than $2n/3$ vertices of $K_n$ with any number of vertices because of the following $2$-coloring: partition the vertices into three almost equal parts $A_1,A_2,A_3$ and color the edges inside $A_i$ and between $A_i,A_{i+1}$ with color $i$. In this coloring no set dominates more than $2n/3$ vertices in any of the three colors. Recently Kral, Liu, Sereni, Whalen and  Yilma \cite{KLSWY} got very close to the solution of Problem \ref{threecol}, showing that $22$ can be replaced by $4$.

\subsection{Covering by monochromatic cycles}

Let the cycle partition number be the minimum $k$ such that the vertex set of any $r$- edge-colored complete graph can be covered by at most $k$ vertex disjoint monochromatic cycles. 

\begin{itemize} \item (\cite{EGYP}) The cycle partition number of any $r$-colored complete graph depends only on $r$ and it is at most $cr^2\log{r}$. \end{itemize}

\begin{conjecture}\pont \label{cycpart} (\cite{EGYP}) The cycle partition number of any $r$-colored complete graph is at most $r$.
\end{conjecture}

The case $r=2$ in Conjecture \ref{cycpart} is due to J. Lehel and was proved for large enough complete graphs by \L uczak, R\"{o}dl and Szemer\'edi in \cite{LRS} using the regularity lemma and later by Allen in \cite{ALL} without it. Then Bessy and Thomass\'e \cite{BT} proved it for all complete graphs with an elementary approach. The estimate $cr^2\log{r}$ of \cite{EGYP} is improved to $cr\log{r}$ in \cite{GRSSz4}. Although Conjecture \ref{cycpart} for $r=3$  was proved asymptotically in \cite{GRSSz3},  Pokrovskiy \cite{PO} found a counterexample. Nevertheless, in the counterexample for $r=3$ three monochromatic cycles cover all but one vertices thus a slightly weaker version of Conjecture \ref{cycpart} can be easily true.  Extensions, related problems can be found in \cite{GYS}.

\subsection{$B+M$ graphs}

At the conference in Orsay in 1976 we had a chat with Paul about $4$-critical ($4$-chromatic but removing any edge becomes $3$-colorable) graphs that can be written as the union of a bipartite graph and a matching \cite{EORS}. We returned to this problem in \cite{CEGS} calling them $4$-critical $B+M$-graphs. A $B+3$ graph is a graph which can be written as the union of a bipartite graph and a matching with three edges.

\begin{problem}\pont (\cite{CEGS}) ``We know that one has to be careful with conjectures in this area. That is why we only suspect that $4$-critical $B+3$-graphs
on $n$ vertices must have at least $2n$ edges asymptotically and dare to conjecture only that they have significantly more than $5n/3$ edges.''
\end{problem}

\section{Szentendre}

During the years 1993 - 1996 Paul spent with us some summer weeks as our guest at Szentendre.  Almost all essentials for his life had been present (a mathematician, a mathematician's wife who could prepare beef Stroganoff and take part in literary and theological debates) - although at the beginning there was no phone and it was a regular program to walk to the nearest phone booth. Another program was to walk to vista point Kada at a hilltop nearby. Or just walk in the garden and enjoy the shade under the huge old walnut tree. Paul often invited us to dinner at the Merry Monks where we became regulars, one waiter have always greeted him asking ``how are you and how are the prime numbers?''

He liked to sit at the terrace in front of the house struggling with problems after problems. Time to time he exclaimed: ``It is very annoying that we do not see this!'' Sometimes we were more successful: ``This is enough for a paper, don't you think so?''

\subsection{Decreasing the diameter of triangle-free graphs}

For a triangle-free graph $G$ we defined $h_d(G)$ (\cite{EGYdiam}) as the minimum number of edges to be added to $G$ to obtain a triangle-free graph of diameter at most $d$.

\begin{itemize}\item (\cite{EGYdiam}) For every connected triangle-free graph $G$ on $n$ vertices, $h_3(G)\le n-1$ and $h_5(G)\le {n-1\over 2}$. \end{itemize}

The case $d=3$ is left open, in particular we asked the following.

\begin{problem}\pont(\cite{EGYdiam}) Is there a positive $\epsilon$ such that $h_4(G)\le (1-\epsilon)n$ for every connected triangle-free graph $G$ on $n$ vertices?
\end{problem}

\subsection{A problem on set mappings}

The first (out of 56) joint paper of Erd\H os and Hajnal, \cite{EH2}, defined set mappings as a function from proper subsets of $S$ to $S$ such that  $f(A)\notin A$. In \cite{EH2} they defined $H(n)$ as the smallest integer for which there exists a set mapping on $S$ with $|S|=n$ such that
$$\cup_{X\subseteq T} f(X)=S$$ for every $T\subset S, |T|\ge H(n)$ and they proved that $\log_2{n}<H(n)$ and conjectured that $H(n)-\log_2(n)$ tends to infinity with $n$. We tried in vain to make the following small step toward this:

\begin{problem}\pont(\cite{EPLAST})  Show that $H(n)>k+1$ for  $n=2^k$.
\end{problem}

\subsection{When every path spans a $3$-colorable subgraph} During the summer of 1995 Paul cited the following from one of his problem books.  ``I asked this with Hajnal: if each odd cycle of a graph spans a subgraph with chromatic number at most $r$ then the chromatic number of the graph is bounded by a function of $r$.'' Some days later I created a warm-up to this question which is still open.

\begin{conjecture}\pont (\cite{GYFRUIT}) If each path of a graph spans at most $3$-chromatic subgraph then the graph is $c$-colorable, perhaps with $c=4$.
\end{conjecture}

\subsection{Balanced colorings} We called an edge coloring of $K_n$ with $r$ colors  {\em balanced} if every subset of $\lfloor n/r \rfloor$ vertices contains at least one edge in each color.
\begin{itemize}\item (\cite{EGYbal}) Balanced $r$-coloring of $K_n$ exists when $n=r^2+r+1$ and $r+1$ is a prime power.
\end{itemize}

We conjectured that this result gives the smallest $n$ for which balanced $r$-coloring exists.

\begin{conjecture}\pont (\cite{EGYbal}) In every $r$-coloring of the edges of $K_{r^2+1}$ there exist $r+1$ vertices with at least one missing color among them ($r\ge 3$, true for $r=3,4$).
\end{conjecture}

\subsection{Nearly bipartite graphs}

Erd\H os and Hajnal asked if every subset $S$ of vertices in a graph $G$ contains an independent set of size at least
$\lfloor {|S|\over 2}\rfloor - k$ then can one remove $f(k)$ vertices from $G$ so that the remaining graph is bipartite?
This is settled in the affirmative by Bruce Reeed \cite{R}. I looked at the case $k=0$ and called a graph $G$ {\em nearly bipartite} if every $S\subseteq V(G)$ has an independent set of size $\lfloor{|S|\over 2}\rfloor$.

\begin{itemize}\item  (\cite{GYFRUIT}) A graph is nearly bipartite if and only if it contains neither two vertex disjoint odd cycles nor odd subdivisions of $K_4$ (where each edge is subdivided to form an odd path).
\end{itemize}

\begin{conjecture}\pont\cite{GYFRUIT}\label{nbip} A nearly bipartite graph can be made bipartite by deleting at most 5 vertices.
\end{conjecture}

Molloy and Reed believed \cite{R} they can prove Conjecture \ref{nbip} but so far they did not work out the details...   The following example shows that Conjecture \ref{nbip} would be sharp. Take a $5\times 5$ grid with vertex set $\{v_{i,j}: 1\le i,j \le 5\}$, subdivide its horizontal edges and add the edges $v_{1,1}v_{5,5}, v_{1,2}v_{5,4},v_{1,3}v_{5,3},v_{1,4}v_{5,2},v_{1,5}v_{5,1}$.

\subsection{Covering by monochromatic paths}

I mentioned to Paul that the vertices of a $2$-colored complete graph can be covered by the vertices of two monochromatic paths - an easy exercise, a footnote in my first paper (with Gerencs\'er) \cite{GGY}. Paul said he did not believe that and it turned out soon that he thought that the covering paths must have the same color. Within two weeks we arrived to a partial answer to this new problem.

\begin{itemize}\item The vertex set of a $2$-colored $K_n$ can be always covered by the vertices of at most $2\sqrt{n}$ monochromatic paths of the same color.
\end{itemize}

\begin{problem}\pont(\cite{EGYCOV})
Is it possible to cover the vertex set of a $2$-colored $K_n$ with at most$\sqrt{n}$ monochromatic paths of the same color? This would be best possible.
\end{problem}

\subsection{$(p,q)$-colorings of complete graphs}

It is sad to look at the submission date on our paper \cite{EGY}:  September  15, 1996 - one day before Paul died in Varsaw.
We called an edge coloring of a complete graph $K_n$ a {\em $(p,q)$-coloring} if every $K_p\subset K_n$ spans at least $q$ colors. Thus a $(p,2)$-coloring means that there is no monochromatic $K_p$, a $(3,3)$-coloring means that the coloring is proper. Let $f(n,p,q)$ denote the minimum number of colors needed for a $(p,q)$-coloring of $K_n$. Some of the general bounds of this paper have been improved by S\'ark\"ozy and Selkow \cite{SSE}. There are much improvements on interesting small cases as well.

On $f(n,4,3)$ the best upper bound is $n^{o(1)}$ given in Mubayi \cite{MU}. The lower bound of  Kostochka and Mubayi \cite{KM} is improved to $c\log{n}$ by Fox and Sudakov \cite{FSU}.  Lower bounds of \cite{EGY} and upper bounds of Mubayi \cite{MU2} imply that  $f(n,4,4)=n^{1/2+o(1)}$.

In the next problem the debate of the authors whether the lower or upper bound is closer to the truth is not resolved (yet). As already noted in \cite{EGY}, $f(n,4,5)\le n-1$ for all even $n\ge 6$ would follow from the existence of a factorization of $K_n$ where there is no $C_4$ in the union of any two color classes. Such factorizations exist if $n=2p$ with odd $p$ or $n=s+1$ where $s$ is not divisible by three.
\begin{problem}\pont(\cite{EGY}) ${5(n-1)\over 6} \le f(n,4,5)\le n$ - improve the estimates!
\end{problem}

One of my favorite problems from \cite{EGY} is to decide whether $f(n,5,9)$ is linear which is equivalent to the following:

\begin{problem}\pont(\cite{EGY}) Is there a constant $c$ such that $K_n$ has a proper edge coloring with $cn$ colors, such that the union of any two color classes has no path or cycle with four edges?
\end{problem}

G\'eza T\'oth proved $2n-6 \le f(n,5,9)$ \cite{Tothg} and  $f(n,5,9)\le 2n^{1+c/\sqrt{log n}}$ is due to Axenovich \cite{AX}.

\subsection{Problems and solutions}
With so many problems to ask, think about, share, transfer, it was unavoidable that sometimes Paul created some confusion.  For me a remarkable adventure of this kind was my problem that Paul announced at the 18-th Southeastern Conference on Combinatorics and Graph Theory held in Boca Raton (February 23-27, 1987).

Years went by, Paul forgot what happened and contributed the problem to himself...  I mentioned his slip of mind at Szentendre in a summer evening of 1995.  ``It is not important who asked the problem, the important thing is that the problem is solved'' - he said and I heartily agreed...


It was not always easy to live and work with Paul. But it was hard to except that we had no more summers with him. Chess, pingpong, card games, conjectures, proofs continued, but not with the same passionate way. And the huge old walnut tree in the garden also died in 1996.
Nevertheless, his results, conjectures, proofs and passionate love of mathematics are with us and carry over to the present and forthcoming generations of mathematicians.


\end{document}